\documentclass[12pt]{article}
\usepackage{amsxtra,amssymb,amsthm,amsmath,latexsym}

\theoremstyle{plain}
\newtheorem{definition}{Definition}[section]
\newtheorem{theorem}{Theorem}[section]

\def\oH{\buildrel\circ\over H}
\def\oH1{\buildrel\circ\over H\kern-.02in{}^1}

\begin{document}

\title{ A simple proof of the Fredholm
alternative and a characterization of the
Fredholm operators.
}

\author{
A.G. Ramm\\
}

\date{}

\maketitle\thispagestyle{empty}

\begin{abstract}
Let $A$ be a linear bounded operator in a Hilbert space $H$,
$N(A)$ and $R(A)$
its null-space and range, and $A^*$ its adjoint. The operator
$A$ is called Fredholm iff
$dim \,N(A)= dim \, N(A^*):=n<\infty$ and $R(A)$ and $R(A^*)$
are closed subspaces of $H$.

A simple  and short proof is given  of the following known
result: $A$ is Fredholm iff $A=B+F$,
where $B$ is an isomorphism and $F$ is a finite-rank operator.
The proof  consists in reduction to a finite-dimensional linear
algebraic system
which is equivalent to the  equation $Au=f$ in the case of
Fredholm operators.

\end{abstract}


\section{Introduction}

{\it The aim of our paper is to prove the
Fredholm alternative and to give a characterization of the class
of Fredholm operators in a very simple way, by a reduction
of the operator equation with a Fredholm operator
to a linear algebraic system in  a finite dimensional space.}

The emphasis is on the simplicity of the
  argument. The paper is written for a wide audience.
The Fredholm alternative is a classical well-known result whose
proof for linear equations of the form $(I+T)u=f$, where $T$ is
a compact operator 
in a Banach space, can be found in 
most of the texts on functional analysis,
of which we mention just  [1]-[2].
 A characterization of the
 set of Fredholm operators one can find in [1], but
it is not given in most of the texts. The proofs in the cited books
follow the classical Riesz 
argument in construction of the Riesz-Fredholm theory.
Though beautiful, this theory is not very simple.

    Our aim is to give  a very short and  simple
 proof of the Fredholm alternative and of a characterization of the class 
of Fredholm operators. This proof is accessible to 
a student with very limited background.
For this reason we give the argument for the case of Hilbert
space, 
but the proof
is  quite easy to adjust for the case of Banach space.

  The idea is to reduce the problem to the one for linear
algebraic systems in finite-dimensional case, 
for which the Fredholm alternative is
a simple fact known to beginners: in a finite-dimensional space
$R^N$ property (1.4) in the Definition 1.1 of Freholm operators
is a consequence of the closedness of any finite-dimensional linear
subspace, since $R(A)$ is such a subspace in $R^N$, while property
(1.3) is a consequence of the simple formulas $r(A)=r(A^*)$ 
and $n(A)=N-r(A),$ valid for matrices, where $r(A)$ is the rank
of $A$ and $n(A)$ is the dimension of the null-space of $A$.

Throughout the paper $A$ is a linear 
bounded operator, $A^*$ its adjoint, $N(A)$ and $R(A)$
are the null-space and the range of $A$.

Recall that an operator $F$ 
with $dim \, R(F)<\infty$ is called  a finite-rank operator,
its rank is equal to  $n:=dim \, R(F)$.

We call a linear bounded operator $B$  on $H$  
an isomorphism if it is
a bicontinuous injection of $H$ onto $H$, that is, $B^{-1}$
is defined on all of $H$ and is bounded.

If $e_j, 1\leq j \leq n$, is an orthonormal basis of  $R(F)$,
then
$Fu=\sum_{j=1}^n(Fu, e_j)e_j$,  so
\begin{equation}
Fu= \sum_{j=1}^n(u, F^* e_j)e_j,
\tag{1.1}\end{equation}
and
\begin{equation}
F^*u= \sum_{j=1}^n(u, e_j) F^* e_j,
\tag{1.2}\end{equation}
where $(u,v)$ is the inner product in $H$.
 
\begin{definition} An operator $A$ is called Fredholm iff ($=$if and only if)
 \begin{equation}
dim\, N(A)= dim \, N(A^*):=n<\infty,
\tag{1.3}
\end{equation}
and
\begin{equation}
R(A)=\overline {R(A)}, \quad R(A^*)=\overline {R(A^*)},
\tag{1.4}\end{equation}
where the overline stands for the closure.
\end{definition}
Recall that
\begin{equation}
H=\overline {R(A)} \oplus N(A^*), \quad  H=\overline {R(A^*)} \oplus N(A) ,
\tag{1.5}\end{equation}
for any linear densely-defined operator $A$, not necessarily
bounded.
For a Fredholm operator $A$ one has:
\begin{equation}
H= {R(A)} \oplus N(A^*), \quad  H= {R(A^*)} \oplus N(A) .
\tag{1.6}
\end{equation}
Consider the equations:
\begin{equation}
Au=f,
\tag{1.7}\end{equation}
\begin{equation}
Au_0=0,
\tag{1.8}\end{equation}
\begin{equation}
A^*v=g,
\tag{1.9}\end{equation}
\begin{equation}
A^*v_0=0.
\tag{1.10}\end{equation}

Let us formulate the Fredholm alternative:
\begin{theorem} If $A=B+F$, where $B$ is an isomorphism and $F$ is
a finite rank operator, then $A$ is Fredholm. 

For any Fredholm operator  $A$ the following (Fredholm)
alternative holds:

1) either (1.8) has only the trivial  solution 
$u_0=0$, and then (1.10)
has only the trivial solution,  and equations (1.7) and (1.9) 
are uniquely solvable for any right-hand sides $f$ and $g$,

or

2) (1.8) has exactly $n>0$ linearly 
independent solutions  $\{\phi_j\}, 1\leq j \leq n$,
and then (1.10) has also $n$ linearly independent solutions 
$\{\psi_j\}, 1\leq j \leq n$, equations (1.7) and (1.9) are solvable iff 
$(f,\psi_j)=0, \, 1 \leq j \leq n$, and correspondingly 
$(g,\phi_j)=0, \, 1 \leq j \leq n$. If they are solvable, 
their solutions are not unique
and their general solutions are respectively:
$u=u_p+ \sum_{j=1}^n a_j \phi_j,$
and 
$v=v_p+ \sum_{j=1}^n b_j \psi_j,$
where $a_j$ and $b_j$ are arbitrary constants, and $u_p$ and $v_p$
are some particular solutions to (1.7) and (1.9), respectively.
\end{theorem}  

Let us give a characterization of the class 
of Fredholm operators, that is,
a necessary and sufficient condition for $A$ to be Fredholm.
\begin{theorem} A linear bounded operator $A$ is Fredholm iff $A=B+F$,
where $B$ is an isomorphism and $F$ has finite rank.
\end{theorem}  

In section 2 we prove these theorems.

\section{Proofs}

\begin{proof}[Proof of Theorem 1.2]
We give a proof of Theorem 1.1 below. From this proof it follows that if
$A=B+F$, where $B$ is an isomorphism and $F$ has finite rank,
then
$A$ is Fredholm.
To prove the converse, 
choose some orthonormal bases $\phi_j$ and $\psi_j, $ in $N(A)$  and 
$N(A^*)$ respectively, using assumption (1.3).
Define
\begin{equation}
Bu:=Au-\sum_{j=1}^n (u,\phi_j) \psi_j:=Au-Fu.
\tag{2.1}\end{equation}
Clearly $F$ has finite rank, and $A=B+F$. Let us prove that $B$
is an isomorphism.
If this is done, then Theorem 1.2 is proved.

We  need to prove that $N(B)=\{0\}$ and $R(B)=H$.
It is known (Banach's theorem), that if $B$ is a linear
injection and $R(B)=H$, then $B^{-1}$ is a bounded operator,
so $B$ is an isomorphism in the sense defined above.

Suppose $Bu=0$. Then $Au=0$  (so that
$u \in N(A)$), and $Fu=0$  (because,
according to (1.6),  $Au$ is orthogonal to $Fu$).
Since $\{\psi_j\},  1 \leq j \leq n$, is a 
linearly independent system,
equation $Fu=0$ implies $(u,\phi_j)=0$ for all $1 \leq j \leq n$,
that is, $u$ is orthogonal to $N(A)$. If $u\in N(A)$ and at the
same
time it is orthogonal to
$N(A)$, then $u=0$.  So, $N(B)={0}$.

Let us now prove that $R(B)=H$:

 Take an arbitrary $f\in H$ and, using 
(1.6), represent it as 
$f=f_1+f_2$ where  $f_1\in R(A)$ and $f_2 \in N(A^*)$
are orthogonal. Thus there is a $u_p\in H$
and some constants  $c_j$
 such that
$f=Au_p + \sum_1^nc_j\psi_j$. We choose 
$u_p$ orthogonal to $N(A)$. This is clearly possible.

We claim that 
$Bu=f,$  where  $u:= u_p - \sum_1^nc_j\phi_j$.
Indeed,   using the orthonormality of the system $\phi_j, \, 1 \leq j \leq n,$
one gets
$Bu=Au_p + \sum_1^nc_j\psi_j=f$.

Thus we have proved that $R(B)=H$.

Theorem 1.2 is proved. 
\end{proof}

We now prove Theorem 1.1.
\begin{proof} [Proof of Theorem 1.1]
 If $A$ is Fredholm, then the statements
 1) and 2) of Theorem 1.1 are equivalent to 
(1.3) and (1.4), since (1.6) follows from (1.4).

Let us prove that if $A=B+F$, where $B$ is an isomorphism and $F$ 
has finite-rank,
then $A$ is Fredholm. Both properties (1.3) and (1.4) are known
for operators in finite-dimensional spaces.
Therefore  we will prove that $A$ is Fredholm if we prove that 
equations (1.7) and (1.9) are equivalent to linear algebraic systems
in a finite-dimensional space.

Let us  prove this equivalence.  We start with equation (1.7),  denote
$Bu:=w$ and get an equivalent to (1.7) equation 

\begin{equation}
w+Tw=f, 
\tag{2.2}\end{equation}
where  $T:=FB^{-1},$ is a finite rank operator  which is of the
same rank $n$ as $F$ 
because $B$ is an isomorphism.
Equation (2.2) is equivalent to (1.7): each solution to (1.7)
is in one-to-one correspondence with a solution of 
(2.2) since $B$ is an isomorphism.
In particular the dimensions of the 
null-spaces $N(A)$ and $N(I+T)$ are equal,
 $R(A)=R(I+T),$  and $R(I+T)$ is closed .
The last claim is a consequence of the  Fredholm alternative for
finite-dimensional
linear equations, but we give an independent  proof
of the closedness of $R(A)$ at the end of the paper.

 Since $T$ is a  finite rank operator,  
the dimension of $N(I+T)$ is finite and  is
not greater than  the rank of $T$. Indeed,
if $u=-Tu$ and $T$ has finite rank $n$, then
$Tu=\sum_{j=1}^n (Tu,e_j)e_j$, where $\{e_j\}_{1\leq j \leq n},$
is an orthonormal basis of $R(T)$, and 
$u=-\sum_{j=1}^n (u,T^*e_j)e_j$, so that
$u$ belongs to a subspace of dimension $n=r(T)$.

Since $A$ and $A^*$ enter symmetrically in 
the statement of Theorem 1.1, 
it is sufficient  to prove  (1.3) and (1.4) for 
$A$ and check that  the dimensions
of $N(A)$ and $N(A^*)$ are equal.

To prove (1.3) and (1.4), let us reduce (1.9)  
to an equivalent equation of the form: 
\begin{equation}
v+T^*v=h, 
\tag{2.3}\end{equation}
where $T^*:=B^{*-1}F^*,$ is the adjoint to $T,$ and
\begin{equation}
h:=B^{*-1}g.
\tag{2.4}\end{equation}
Since $B$ is an isomorphism,  $(B^{-1})^*=(B^*)^{-1}$. 
Applying $B^{*-1}$ to equation (1.9), one gets 
an equivalent equation (2.3) and 
$T^*$  is a finite-rank operator of the same rank  $n$ as $T$.

The last claim is easy to prove: if $\{e_j\}_{1\leq j \leq n}$
is a basis in $R(T),$ then $Tu=\sum_{j=1}^n (Tu,e_j)e_j$,
and $T^*u=\sum_{j=1}^n (u,e_j)T^*e_j$, so 
$r(T^*)\leq r(T)$. By symmetry one has $r(T)\leq r(T^*),$
and the claim is proved.

  Writing explicitly the  linear  algebraic systems, equivalent to the
equations (2.2) and (2.3), one sees that the matrices of these systems
are adjoint. The system equivalent to equation (2.2) is:
\begin{equation}
c_i+\sum_1^n  t_{ij} c_j=f_i, 
\tag{2.5}\end{equation}
where
$$
t_{ij}:=(e_j, T^*e_i), \,\, c_j:=(w,T^*e_j),\,\,
f_i:=(f,T^*e_i),
$$
and the one equivalent to (2.3) is:
\begin{equation}
\xi_i+\sum_1^n  t_{ij}^{*} \xi_j=h_i, 
\tag{2.6}
\end{equation}
where
$$
 t_{ij}^{*}=(T^*e_j, e_i), \,\,
\xi_j:=(v,e_j),\,\, h_i:=(h,e_i),
$$
and $t_{ij}^{*}$ is the matrix adjoint to $t_{ij}$. For
linear algebraic systems  (2.5) and (2.6)  
the Fredholm alternative is a
well-known elementary result. These systems are equivalent to
equations (1.7) and (1.9), respectively. 
Therefore the Fredholm alternative holds for
equations (1.7) and (1.9), 
so that  properties (1.3) and (1.4) are proved.
Theorem 1.1 is proved.
\end{proof}

In conclusion let us explain in detail why, for example,
equations (2.2) and (2.5) are equivalent
in the following sense: every solution to
(2.2) generates a solution to (2.5) and vice versa.

It is clear that (2.2) implies (2.5): just take the inner
product of 
(2.2) with $T^*e_j$ and get (2.5). So, each solution to
(2.2) generates a solution to (2.5). We claim that
each solution to (2.5) generates a solution to (2.2).
Indeed, let
$c_j$ solve (2.5). Define $w:=f-\sum_1^n c_je_j$. 
Then $Tw=Tf-\sum_{j=1}^n c_jTe_j=\sum_{i=1}^n
[(Tf,e_i)e_i- 
\sum_{j=1}^n c_j(Te_j,e_i)e_i]=\sum_{i=1}^n c_ie_i=f-w.$
Here we have used (2.5) and took into account
that $(Tf,e_i)=f_i$ and $(Te_j,e_i)=t_{ij}.$
Thus, the element $w:=f-\sum_1^n c_je_j$ solves (2.2),
as claimed.

It is easy to check that if $\{w_1,....w_k\}$ are $k$
linearly independent solutions to the homogeneous
version of equation (2.2), then the corresponding $k$
solutions $\{c_{1m},.....c_{nm}\}_{1\leq m \leq k}$
of the homogeneous version of the system (2.5)
are also linearly indepenedent, and vice versa.

Let us give an independent proof of property (1.4): 

{\it $R(A)$ is closed
if $A=B+F$, where $B$ is an isomorphism
and $F$ is  a finite rank operator.}

Since  $A=(I+T)B$ and $B$ is an isomorphism, it is sufficient to
prove that $R(I+T)$ is closed if $T$ has finite rank.

Let  $u_j +Tu_j :=f_j\to f$ as $j \to \infty.$
Without loss of generality choose $u_j$ orthogonal to $N(I+T)$.
  We want to prove that there exists
 a $u$ such that 
$(I+T)u=f$.  Suppose first that $\sup _j ||u_j||<\infty.$
Since $T$ is a finite-rank operator, $Tu_j$ converges in $H$
for some subsequence which
is denoted $u_j$ again.  (Recall that in finite-dimensional spaces bounded sets are
precompact).
This implies that $u_j=f_j-Tu_j $ converges in $H$ to an element $u$.
Passing to the limit, one gets $(I+T)u=f.$
 To complete the proof, let us establish that 
$\sup _j ||u_j||<\infty$ . Assuming   that this is
false, one can choose a  subsequence, denoted $u_j$ again, 
such that  $||u_j||>j.$
Let  $z_j:=u_j/||u_j||$. Then  $||z_j||=1,$  
$z_j$ is orthogonal to $N(I+T),$ and 
$z_j+Tz_j=f_j/||u_j|| \to 0$.
As above, it follows that $z_j \to z$ in $H$,  
and passing to the limit
in the equation for $z_j$ one gets $z+Tz=0$. 
Since $z$ is orthogonal to $N(I+T)$,
it follows that $z=0$. This is a contradiction 
since $||z||=lim_{j\to \infty} ||z_j||=1$.
This contradiction proves the desired estimate 
and the proof is completed.

The above proof is valid for any compact linear operator $T$. 
If $T$ is a finite-rank operator,
then the closedness of  $R(I+T)$ follows 
from a simple observation: finite-dimensional linear 
spaces are closed.

\section*{Acknowledgement}
The author thanks Professor R.Burckel for very useful comments.
\newpage

\noindent {\Large{\bf References}}

\vspace{.1in}

\noindent 1.~Kantorovich, L., Akilov, G.,  Functional
analysis in normed spaces,
Macmillan, New York, 1964

\vspace{.1in}

\noindent 2.~Rudin, W.,  Functional analysis,
McGraw Hill, New York, 1973

$$
$$

\noindent Mathematics Department, Kansas State University, \\
 Manhattan, KS 66506-2602, USA\\
 ramm@math.ksu.edu\\
http://www.math.ksu.edu/\,$\widetilde{\ }$\,ramm

\end{document}